\documentclass[12pt]{article}
\usepackage{amssymb,amsmath,amsthm,amsfonts,enumerate, bbm, mathdots}

\usepackage{latexsym}
\usepackage{amscd}
\usepackage{stmaryrd}
\usepackage{color}
\usepackage[all]{xypic}
\usepackage{epsfig}
\usepackage{graphics}
\usepackage{ifthen}
\usepackage{varioref}
\usepackage{rotating}
\usepackage{extarrows}
\usepackage{cite}
\usepackage{mathrsfs}

\numberwithin{equation}{section}

\theoremstyle{plain}
\newtheorem{thm}{Theorem}[section]
\newtheorem{cor}[thm]{Corollary}
\newtheorem{lem}[thm]{Lemma}
\newtheorem{prop}[thm]{Proposition}
\newtheorem{rem}[thm]{Remark}
\newtheorem{Def}[thm]{Definition}
\newtheorem{eg}[thm]{Example}

\definecolor{darkgreen}{rgb}{0.0625,0.64,0.0625}
\usepackage[
            pdfstartview=FitH,
            CJKbookmarks=true,
            bookmarksnumbered=true,
            bookmarksopen=true,
            colorlinks,
            pdfborder=001,
            linkcolor=blue,
            anchorcolor=green,
            citecolor=red
            ]{hyperref}

\newfont{\scyr}{wncyr10 scaled 550}

\def\proof{\noindent {\bf Proof.\;}}

\def\id{\operatorname{id}}
\def\Spec{\operatorname{Spec}}

\allowdisplaybreaks

\begin{document}

\title{Rota-Baxter $C^{\ast}$-algebras}

\date{\small ~ \qquad\qquad School of Mathematical Sciences, Tongji University \newline No. 1239 Siping Road,
Shanghai 200092, China}

\author{Zhonghua Li\thanks{E-mail address: zhonghua\_li@tongji.edu.cn} ~and ~Shukun Wang\thanks{E-mail address: 2010165@tongji.edu.cn}}

\maketitle

\begin{abstract}
This paper introduces the notion of Rota-Baxter $C^{\ast}$-algebras. Here a Rota-Baxter $C^{\ast}$-algebra is a $C^{\ast}$-algebra with a Rota-Baxter operator. Symmetric Rota-Baxter operators, as special cases of Rota-Baxter operators on $C^{\ast}$-algebra, are defined and studied. A theorem of Rota-Baxter operators on concrete $C^{\ast}$-algebras is given, deriving the relationship between two kinds of Rota-Baxter algebras. As a corollary, some connection between $\ast$-representations and Rota-Baxter operators is given. The notion of representations of Rota-Baxter $C^{\ast}$-algebras are constructed, and a theorem of representations of direct sums of Rota-Baxter representations is derived. Finally using Rota-Baxter operators,  the notion of quasidiagonal operators on $C^{\ast}$-algebra is reconstructed.
\end{abstract}

{\small

}


\section{Introduction}\label{Sec:Intro}
\newcommand{\RNum}[1]{\uppercase\expandafter{romannumeral #1\relex}}

Let $\mathcal{A}$ be an associative algebra over a given field $F$. A linear operator $P$ on $\mathcal{A}$ is called a Rota-Baxter operator of weight $\lambda\in F$ if it satisfies
\begin{equation}\label{1}
	P(a)P(b)=P(aP(b))+P(P(a)b)+\lambda P(ab),\qquad\forall a,\ b\in \mathcal{A}.
\end{equation}
An associative algebra with a Rota-Baxter operator is called a Rota-Baxter algebra, which can be regarded as an analogue of a differential algebra. In fact, when taking $\lambda=0$, Eq. \eqref{1} is an algebraic abstraction of the formula of integration by parts.

The study of Rota-Baxter algebras originated from probability theory \cite{Baxter1} and has  found applications in many areas of mathematics and physics. In the late 1990s, Rota-Baxter algebras were found their applications in the work of Connes and Kremier \cite{Connes1} regarding the renormalization of perturbative quantum field theory. Since 2000, the connection between the classical Yang-Baxter equation in mathematical physics and Rota-Baxter operators has been found in \cite{Augiar}.

The representations of Rota-Baxter algebras were studied in \cite{Linzongzhu1}, where some basic concepts and properties were established. However, it is still in the early stages of development. In \cite{Linzongzhu1}, Lin and Qiao studied the representations and regular-singular decomposition of Laurent series Rota-Baxter algebras. In \cite{Liguo3}, regular-singular decomposition of Rota-Baxter modules were obtained under the condition of quasi-idempotency. And in \cite{Liqiao1}, representations of the polynomial Rota-Baxter algebras were studied.

The theory of operator algebras in Hilbert spaces was initiated by von Neumann \cite{von} in 1929. In \cite{Murray}, Murray and von Neumann laid the foundation of the theory of $W^{\ast}$-algebras. The notion of $C^{\ast}$-algebras was introduced by Gelfand and Naimark \cite{Neumark} in 1943. Basic theory of representations of $C^{\ast}$-algebras was established in \cite{Gelfand, Kadison, Rtprosser, Segal}. Derivations on $C^{\ast}$-algebras and $W^{\ast}$-algebras were studied deeply in the 1960s and 1970s in \cite{Kadison,Rickart,Sakai,Kadison1}. Especially, in \cite{Sakai}, it was shown that all derivations on $W^{\ast}$-algebras are inner derivations. While Rota-Baxter operators on $C^{\ast}$-algebras have not been studied.

In this note, we establish basic concepts of Rota-Baxter $C^{\ast}$-algebras. Some special Rota-Baxter operators in $C^{\ast}$-algebras, which are called symmetric Rota-Baxter operators and Rota-Baxter operators matching projections are studied. And representations of Rota-Baxter $C^{\ast}$-algebras are also established, which are special cases of traditional Rota-Baxter representations.

A $C^{\ast}$-algebra $\mathcal{A}$ is a special algebra in the complex field. Combined with a Rota-Baxter operator on $\mathcal{A}$, the basic concepts of Rota-Baxter $C^{\ast}$-algebras were established in Subsection \ref{2.1}. A Rota-Baxter operator $P$ on the $C^{\ast}$-algebra $\mathcal{A}$ is symmetric if $P(a^{\ast})=P(a)^{\ast}$ for any $a\in\mathcal{A}$. In Subsection \ref{2.2}, We show that a $C^{\ast}$-algebra $\mathcal{A}$ can be decomposed into a direct sum of two $C^{\ast}$-subalgebras if and only if there is a bounded idempotent symmetric Rota-Baxter operator of weight $-1$ on $\mathcal{A}$.

In Section \ref{3}, we study Rota-Baxter operators on the $C^{\ast}$-algebra $B(\mathcal{H})$, where $\mathcal{H}$ is a Hilbert space.  In Subsection \ref{3.1}, the notion of Rota Baxter operators that match projections are introduced. We find that the Rota-Baxter operators matching projections on Hilbert spaces have many good properties. In fact, we can construct Rota-Baxter operators of this kind from Rota-Baxter operators on $C^{\ast}$-subalgebras of $B(\mathcal{H})$. As a corollary, we study the relationship between the Rota-Baxter operators matching projections and $\ast$-representations of $C^{\ast}$-algebras on a Hilbert space.  In Subsection \ref{3.2}, we construct representations of Rota-Baxter $C^{\ast}$-algebras.

At last in Section \ref{4}, we reconstruct the notion of quasidiagonal operators of $C^{\ast}$-algebras with the help of Rota-Baxter operators.


\section{Rota-Baxter $C^{\ast}$-algebras}

\subsection{$C^{\ast}$-algebras}\label{2.0}

We first recall some basic concepts of $C^{\ast}$-algebras from \cite{Sakai1}. Let $F$ be the complex field $\mathbb{C}$ or the real number field $\mathbb{R}$.

Let $\mathcal{A}$ be an associative algebra over $F$. The algebra $\mathcal{A}$ is called a normed algebra if associated to each element $a$ in $\mathcal{A}$ there is a real number $\|a\|$, called the norm of $a$, with the properties:
\begin{description}
\item[(i)] $\|a\|\ge 0$ for any $a\in\mathcal{A}$, and $\|a\|=0$ if and only if $a=0$;
\item[(ii)] $\|a+b\|\le\|a\|+\|b\|$ for any $a,b\in\mathcal{A}$;
\item[(iii)] $\|\lambda a\|=|\lambda|\|a\|$ for any $a\in\mathcal{A}$ and $\lambda\in F$;
\item[(iv)] $\|ab\|\le\|a\|\|b\|$ for any $a,b\in\mathcal{A}$.
\end{description}
The topology defined by the norm $\|\cdot\|$ on $\mathcal{A}$ is called the uniform topology. If $\mathcal{A}$ is complete with respect to the norm, then $\mathcal{A}$ is called a Banach algebra. A map $\mathcal{A}\rightarrow \mathcal{A}; a\mapsto a^{\ast}$ is called an involution if it satisfies the following conditions:
\begin{description}
\item[(i)] $(a^{\ast})^{\ast}=a$ for any $a\in\mathcal{A}$;
\item[(ii)] $(a+b)^{\ast}=a^{\ast}+b^{\ast}$ for any $a,b\in\mathcal{A}$;
\item[(iii)] $(\lambda a)^{\ast}=\overline{\lambda}a^{\ast}$ for any $a\in\mathcal{A}$ and $\lambda\in F$;
\item[(iv)] $(ab)^{\ast}=b^{\ast}a^{\ast}$ for any $a,b\in\mathcal{A}$.
\end{description}
An algebra with an involution $\ast$ is called a $\ast$-algebra. Finally, a Banach $\ast$-algebra $\mathcal{A}$ is called a ($F$-linear) $C^{\ast}$-algebra if it satisfies
$$\|a^{\ast}a\|=\|a\|^{2}$$
for any $a\in\mathcal{A}$. Note that a $\mathbb{C}$-linear $C^{\ast}$-algebra is natural a $\mathbb{R}$-linear $C^{\ast}$-algebra.

A subset $S$ of a $C^{\ast}$-algebra $\mathcal{A}$ is called self-adjoint if $a^{\ast}\in S$ for any $a\in S$. In particular, an element $a\in\mathcal{A}$ is called self-adjoint if $a^{\ast}=a$. A self-adjoint, uniformly closed subalgebra of $\mathcal{A}$ is called a $C^{\ast}$-subalgebra of $\mathcal{A}$, which is also a $C^{\ast}$-algebra.

As there is a new structure $\ast$ on a $C^{\ast}$-algebra, it is natural to consider the $\ast$-homomorphisms between two $C^{\ast}$-algebras. Let $\mathcal{A}$ and $\mathcal{B}$ be two $C^{\ast}$-algebras. A map $\phi:\mathcal{A}\rightarrow \mathcal{B}$ is called a $\ast$-homomorphism if it satisfies
\begin{description}
\item[(i)] $\phi(a+b)=\phi(a)+\phi(b)$ for any $a,b\in\mathcal{A}$;
\item[(ii)] $\phi(\lambda a)=\lambda\phi(a)$ for any $a\in\mathcal{A}$ and $\lambda\in F$;
\item[(iii)] $\phi(ab)=\phi(a)\phi(b)$ for any $a,b\in\mathcal{A}$;
\item[(iv)] $\phi(a^{\ast})=\phi(a)^{\ast}$ for any $a\in\mathcal{A}$.
\end{description}

We recall the definition of the direct sum of $C^{\ast}$-algebras. Let $\{\mathcal{A}_{k}\}_{k\in \Lambda}$ be a family of $C^{\ast}$-algebras. We define the direct sum
$$\bigoplus\limits_{k\in \Lambda}\mathcal{A}_{k}=\left\{(a_{k})_{k\in \Lambda}\mid a_{k}\in\mathcal{A}_{k} \;\text{for any\;}k\in\Lambda, \underset{k\in \Lambda}{\sup}\|a_{k}\|<\infty\right\}.$$
Then $\bigoplus\limits_{k\in \Lambda}\mathcal{A}_{k}$ is a $C^{\ast}$-algebra under the following operators:
\begin{description}
\item[(i)] $(a_{k})+(b_{k})=(a_{k}+b_{k})$;
\item[(ii)] $\lambda(a_{k})=(\lambda a_{k}),\qquad(\lambda\in F)$;
\item[(iii)] $(a_{k})(b_{k})=(a_{k}b_{k})$;
\item[(iv)] $\|(a_{k})\|=\underset{k}{\sup}\|a_{k}\|$;
\item[(v)] $(a_{k})^{\ast}=(a_{k}^{\ast})$.
\end{description}
For any $k\in\Lambda$, it is natural to regard $\mathcal{A}_{k}$ as a $C^{\ast}$-subalgebra of $\bigoplus\limits_{k\in \Lambda}\mathcal{A}_{k}$.

\subsection{Rota-Baxter $C^{\ast}$-algebras}\label{2.1}

Now we introduce the concept of Rota-Baxter $C^{\ast}$-algebras.

\begin{Def}\label{def crb}
Let  $\mathcal{A}$ be a $C^{\ast}$-algebra. A linear operator $P:\mathcal{A}\to\mathcal{A}$ is called a Rota-Baxter operator of weight $\lambda\in F$ on  $\mathcal{A}$  if it satisfies:
\begin{equation*}
	P(a)P(b)=P(aP(b))+P(P(a)b)+\lambda P(ab), \qquad\forall a,\ b\in \mathcal{A}.
\end{equation*}
\end{Def}

If $P$ is a Rota-Baxter operator of weight $\lambda$, then it is easy to verify that
$$\widetilde{P}=-\lambda \id_{\mathcal{A}}-P$$
is also a Rota-Baxter operator of weight $\lambda$.

\begin{Def}\label{crbo}
A Rota-Baxter $C^{\ast}$-algebra of weight $\lambda$ is a pair $(\mathcal{A},P)$ consisting of a $C^{\ast}$-algebra $\mathcal{A}$, and a Rota-Baxter operator $P:\mathcal{A}\to\mathcal{A}$ of weight $\lambda$ on $\mathcal{A}$.
\end{Def}

We give an example of Rota-Baxter $C^{\ast}$-algebras.

\begin{eg}\label{eg1}
Let $C([0,1])$ be the set of all complex valued, continuous functions on the closed interval $[0,1]$. Then $C([0,1])$ becomes a commutative algebra over $\mathbb{C}$ under pointwise addition and multiplication.  It was showed in \cite{Sakai1} that $C([0,1])$ is a commutative $C^{\ast}$-algebra, where for any $f\in C([0,1])$, the norm of $f$ is defined by
$$\|f\|=\underset{x\in [0,1]}{\sup}|f(x)|,$$
and the involution $f^\ast$ of $f$ is given by
$$f^{\ast}(x)=\overline{f(x)},\quad \forall x\in [0,1].$$
The linear operator $T$ on $C([0,1])$ defined by
$$T(f)(x)=\int_{0}^{x}f(s)ds,\quad \forall f\in C([0,1]),\;\forall x\in[0,1]$$
is called the Volterra operator. It is easy to verify that $T$ is a Rota-Baxter operator of weight $0$ on $C([0,1])$.
Then $(C([0,1]),T)$ is a Rota-Baxter $C^{\ast}$-algebra.
\end{eg}

As in \cite{Sakai1}, basic concepts on $C^{\ast}$-algebras can be similarly defined for Rota-Baxter $C^{\ast}$-algebras. Particularly, a Rota-Baxter $C^{\ast}$-subalgebra of a Rota-Baxter $C^{\ast}$-algebra $\mathcal{A}$ is a $C^{\ast}$-subalgebra $I$ of $\mathcal{A}$ such that $P(I)\subseteq I$. A Rota-Baxter $C^{\ast}$-algebra homomorphism $\phi:(\mathcal{A}_{1},P_{1})\to(\mathcal{A}_{2},P_{2})$ between two Rota-Baxter $C^{\ast}$-algebras of the same weight $\lambda$ is a $\ast$-homomorphism such that $\phi\circ P_{1}=P_{2}\circ\phi$.

In below, we always take $F=\mathbb{C}$. Hence all the Hilbert spaces and algebras are assumed to be over the complex field $\mathbb{C}$.

We recall the definition of derivations of $C^{\ast}$-algebras from \cite{Sakai1}. Let $\mathcal{A}$ be a $C^{\ast}$-algebra. A linear map $\delta:\mathcal{A}\rightarrow \mathcal{A}$ is called a derivation if $$\delta(ab)=\delta(a)b+a\delta(b),\qquad\forall a,\ b\in\mathcal{A}.$$
It is easy to see that if $P$ is an invertible Rota-Baxter operator of weight $0$ on $\mathcal{A}$, then the inverse operator of $P$ is a derivation of $\mathcal{A}$.
As an application, we have the following proposition.

\begin{prop}
If $\mathcal{A}$ is a commutative $C^{\ast}$-algebra, then there is no invertible Rota-Baxter operator of weight $0$ on $\mathcal{A}$.
\end{prop}

\proof
From \cite{Sakai1}, we know that if $\delta$ is a derivation on $\mathcal{A}$, then $\delta=0$.
 While if  there is an invertible Rota-Baxter operator $P$ of weight $0$ on $\mathcal{A}$, then the inverse of $P$ is an invertible derivation on $\mathcal{A}$, which is a contradiction.
\qed

\subsection{ Symmetric Rota-Baxter operators}\label{2.2}

In this subsection, we introduce the notion of symmetric Rota-Baxter operators.
\begin{Def}
	Let $\mathcal{A}$ be a $C^{\ast}$-algebra and let $\mathcal{A}_1$ be a $C^{\ast}$-subalgebra of $\mathcal{A}$. A Rota-Baxter operator $P$ on $\mathcal{A}$ is called symmetric on $\mathcal{A}_1$ if
$$P(a^{\ast})=P(a)^{\ast}$$
for any $a\in\mathcal{A}_1$. If $P$ is symmetric on $\mathcal{A}$, we just call that $P$ is symmetric.
\end{Def}

\begin{eg}
The Volterra operator $T$ on $C([0,1])$ in Example \ref{eg1} is symmetric.
\end{eg}

Let $\mathcal{A}$ be a commutative $C^{\ast}$-algebra. Let $\mathcal{A}^{s}$ be the set of all self-adjoint elements in $\mathcal{A}$. It is easy to verify that $\mathcal{A}^{s}$ is a $\mathbb{R}$-linear subalgebra of $\mathcal{A}$. We have the following result.

\begin{prop}
Let $\mathcal{A}$ be a commutative $C^{\ast}$-algebra and $\lambda\in \mathbb{R}$. Then there is a $1$-$1$ correspondence between the following two sets:
\begin{description}
\item[(1)] the set of symmetric Rota-Baxter operators of weight $\lambda$ on  $\mathcal{A}$;
\item[(2)] the set of $\mathbb{R}$-linear Rota-Baxter operators of weight $\lambda$ on $\mathcal{A}^{s}$.
\end{description}
\end{prop}

\proof
Denote by $\Gamma$ the set of symmetric Rota-Baxter operators of weight $\lambda$ on  $\mathcal{A}$, and by $\Gamma_1$ the set of $\mathbb{R}$-linear Rota-Baxter operators of weight $\lambda$ on $\mathcal{A}^{s}$.

Let $P\in \Gamma$. Then for any $a\in \mathcal{A}^{s}$, we have
$$P(a)=P(a^{\ast})=P(a)^{\ast}.$$
Hence we get a restriction $P|_{\mathcal{A}^{s}}:\mathcal{A}^{s}\rightarrow \mathcal{A}^{s}$, which is an element of $\Gamma_1$. Therefore we have a map $\Phi:\Gamma\rightarrow\Gamma_1$, such that $\Phi(P)=P|_{\mathcal{A}^{s}}$.

Conversely, let $P_1\in\Gamma_1$. For any $a\in\mathcal{A}$, we have $a=a_1+i a_2$, where $i=\sqrt{-1}$ and
$$a_{1}=\frac{1}{2}(a+a^{\ast})\in\mathcal{A}^s,\qquad a_{2}=\frac{1}{2i}(a-a^{\ast})\in\mathcal{A}^s.$$
Then we get a map $P:\mathcal{A}\rightarrow\mathcal{A}$ by setting $P(a)=P_1(a_1)+iP_1(a_2)$. It is obvious that $P$ is linear.

We show that  $P$ satisfies the Rota-Baxter relation. Let $b=b_1+ib_2$ be another element of $\mathcal{A}$ with $b_1,b_2\in\mathcal{A}^s$ defined similarly as above. We have
\begin{align*}
P(a)P(b)=&(P_1(a_{1})+iP_1(a_{2}))(P_1(b_{1})+iP_1(b_{2}))\\
=&P_1(a_{1})P_1(b_{1})-P_1(a_{2})P_1(b_{2})+i(P_1(a_{2})P_1(b_{1})+P_1(a_{1})P_1(b_{2})).
\end{align*}
Since $P_1$ is a Rota-Baxter operator on $\mathcal{A}^s$, we find
\begin{align*}
P(a)P(b)
=&P_1(a_{1}P_1(b_{1})+P_1(a_{1})b_{1}+\lambda a_{1}b_{1})-P_1(a_{2}P_1(b_{2})+P_1(a_{2})b_{2}+\lambda a_{2}b_{2})\\
&+iP_1(a_{2}P_1(b_{1})+P_1(a_{2})b_{1}+\lambda a_{2}b_{1})+iP_1(a_{1}P_1(b_{2})+P_1(a_{1})b_{2}+\lambda a_{1}b_{2})\\
=&P_1(a_{1}P_1(b_{1}))-P_1(a_{2}P_1(b_{2}))+iP_1(a_{2}P_1(b_{1}))+iP_1(a_{1}P_1(b_{2}))\\
&+P_1(P_1(a_{1})b_{1})-P_1(P_1(a_{2})b_{2})+iP_1(P_1(a_{2})b_{1})+iP_1(P_1(a_{1})b_{2})\\
&+\lambda P_1(a_{1}b_{1}- a_{2}b_{2})+i\lambda P_1(a_{2}b_{1}+a_{1}b_{2}).
\end{align*}
According to the following decompositions
\begin{align*}
&ab=(a_1b_1-a_2b_2)+i(a_2b_1+a_1b_2),\\
&aP(b)=(a_1P_1(b_1)-a_2P_1(b_2))+i(a_2P_1(b_1)+a_1P_1(b_2)),\\
&P(a)b=(P_1(a_1)b_1-P_1(a_2)b_2)+i(P_1(a_2)b_1+P_1(a_1)b_2),
\end{align*}
we finally get
$$P(a)P(b)=P(aP(b))+P(P(a)b)+\lambda P(ab),$$
which means that $P$ is a Rota-Baxter operator of weight $\lambda$ on  $\mathcal{A}$.

Now we show that $P$ is symmetric.  Since we have the decomposition $a^{\ast}=a_1-ia_2$, we get
$$P(a^{\ast})=P_1(a_1)-iP_1(a_2).$$
On the other hand, we have
$$P(a)^{\ast}=(P_1(a_1)+iP_1(a_2))^{\ast}=P_1(a_1)^{\ast}-iP_1(a_2)^{\ast}=P_1(a_1)-iP_1(a_2).$$
Hence $P(a^{\ast})=P(a)^{\ast}$ and then $P$ is symmetric.

Therefore we get a map $\Psi:\Gamma_1\rightarrow\Gamma$, such that $\Psi(P_1)=P$.

Finally, it is easy to verify that $\Phi\circ\Psi$ and $\Psi\circ\Phi$ are identity maps. And the proposition is proved.
\qed

The following result gives an equivalent condition for the idempotent and symmetric Rota-Baxter operators of weight $-1$ on $C^{\ast}$-algebras.

\begin{thm}
Let $\mathcal{A}$ be a $C^{\ast}$-algebra and let $P$ be a bounded linear operator on $\mathcal{A}$. Then the following two statements are equivalent:
\begin{description}
\item[(1)]$P$ is an idempotent and symmetric Rota-Baxter operator of weight $-1$ on $\mathcal{A}$;
\item[(2)] There is a direct sum decomposition $\mathcal{A}=\mathcal{A}_{1}\oplus\mathcal{A}_{2}$, such that $\mathcal{A}_{1}$ and $\mathcal{A}_{2}$ are $C^{\ast}$-subalgebras of $\mathcal{A}$ and $P$ is the projection of $\mathcal{A}$ onto $\mathcal{A}_{1}$.
\end{description}
\end{thm}

\proof
Assume that $P$ is an idempotent and symmetric Rota-Baxter operator of weight $-1$ on $\mathcal{A}$. Let $\mathcal{A}_1=P(\mathcal{A})$ and $\mathcal{A}_2=\widetilde{P}(\mathcal{A})$, where $\widetilde{P}=\id_{\mathcal{A}}-P$. Then from\cite[Theorem 1.1.13]{L.G1}, we know that $\mathcal{A}_1$ and $\mathcal{A}_2$ are subalgebras of $\mathcal{A}$, $\mathcal{A}=\mathcal{A}_{1}\oplus\mathcal{A}_{2}$ is a direct sum decomposition, and $P$ is the projection of $\mathcal{A}$ onto $\mathcal{A}_1$. For any $a=P(b) \in \mathcal{A}_1$ with $b\in\mathcal{A}$, we have
$$a^{\ast}=P(b)^{\ast}=P(b^{\ast})\in \mathcal{A}_1.$$
Hence $\mathcal{A}_1$ is closed under the involution $\ast$ and then self-adjoint. Now let $\{P(b_{n})\}$ be a sequence in $\mathcal{A}_1$ that uniformly converges to $a\in \mathcal{A}$. Then since $P$ is bounded and idempotent, we get $\{P(b_n)\}$ uniformly converges to $P(a)\in\mathcal{A}_1$. Therefore $\mathcal{A}_1$ is a $C^{\ast}$-subalgebra of $\mathcal{A}$. It is similar to show that $\mathcal{A}_2$ is a $C^{\ast}$-subalgebra of $\mathcal{A}$.

Conversely, assume that we have a direct sum decomposition $\mathcal{A}=\mathcal{A}_1\oplus\mathcal{A}_2$, such that $\mathcal{A}_1$ and $\mathcal{A}_2$ are $C^{\ast}$-subalgebras of $\mathcal{A}$, and $P$ is the projection of $\mathcal{A}$ onto $\mathcal{A}_{1}$. From \cite[Theorem 1.1.13]{L.G1}, we have $P$ is an idempotent Rota-Baxter operator of weight $-1$ on $\mathcal{A}$. For any $a\in\mathcal{A}$, assume that $a=a_1+a_2$ with $a_1\in\mathcal{A}_1$ and $a_2\in\mathcal{A}_2$, then we have $a^{\ast}=a_1^{\ast}+a_2^{\ast}$. Since $\mathcal{A}_i$ is a $C^{\ast}$-subalgebra of $\mathcal{A}$, we get $a_i^{\ast}\in\mathcal{A}_i$. Hence
$$P(a^{\ast})=a_1^{\ast}=P(a)^{\ast},$$
which means that $P$ is symmetric. Finally, from \cite[Corollary 1.2.6]{Sakai1} and the definition of the direct sum of $C^{\ast}$-algebras in Subsection \ref{2.0}, we know that
$$\|a\|=\max(\|a_{1}\|,\|a_{2}\|)\ge \|a_1\|=\|P(a)\|.$$
Hence we have
$$\|P\|=\underset{a\in\mathcal{A},a\neq 0}{\sup}\frac{\|P(a)\|}{\|a\|}\le 1,$$
which means that $P$ is bounded.
\qed


\section{Rota-Baxter operators matching projections on Hilbert spaces}\label{3}

\subsection{The $C^\ast$-algebra $B(\mathcal{H})$}

We recall the notion of Hilbert spaces. A complex linear space $\mathcal{H}$ is an inner product space if associated to any $x,y\in\mathcal{H}$ there is a product $\left<x,y\right>$  which satisfies the following properties:
\begin{description}
	\item[(i)] $\left<x,y\right>=\overline{\left<y,x\right>}$ for any $x,y\in\mathcal{H}$;
	\item[(ii)] $\left<\lambda_1x_{1}+\lambda_2x_{2},y\right>=\lambda_{1}\left<x_{1},y\right>+\lambda_{2}\left<x_{2},y\right>$ for any $\lambda_{1},\lambda_{2}\in\mathbb{C}$ and $x_1,x_2,y\in\mathcal{H}$;
	\item[(iii)] $\left<x,x\right>=0$ for $x=0$ and $\left<x,x\right>> 0$ for any nonzero $x$ in $\mathcal{H}$.
\end{description}
Then the norm of $\mathcal{H}$ induced from the inner product is defined by $\|x\|=\sqrt{\left<x,x\right>}$ for any $x\in\mathcal{H}$. We call $\mathcal{H}$ a Hilbert space if it is complete with respect to the norm.

Let $B(\mathcal{H})$ be the set of all bounded linear operators on the Hilbert space $\mathcal{H}$. For any $f\in B(\mathcal{H})$ and any $x\in\mathcal{H}$, we sometimes denote $f(x)$ by $f.x$. From \cite{Sakai1}, we know that $B(\mathcal{H})$ is a $C^{\ast}$-algebra.  In fact, the algebra $B(\mathcal{H})$ is a Banach algebra with the norm given by
$$\|f\|=\underset{\|x\|\le1, x\in\mathcal{H}}{\sup}\|f(x)\|,\quad \forall f\in B(\mathcal{H}).$$
And then $B(\mathcal{H})$ is a $C^{\ast}$-algebra with the involution $\ast$ determined  by
$$\left<f(x),y\right>=\left<x,f^{\ast}(y)\right>,\quad \forall f\in B(\mathcal{H}), \quad \forall x,y\in\mathcal{H}.$$

For a projection $p$ on the Hilbert space $\mathcal{H}$, we mean that $p$ is a linear operator on $\mathcal{H}$ satisfying $p^2=p$ and $p^{\ast}=p$. We set $p^{\perp}=\id_{H}-p$. Then we have $\mathcal{H}=\mathcal{H}_{1}\oplus\mathcal {H}_{2}$, where $\mathcal{H}_{1}=p\mathcal{H}$ and $\mathcal{H}_{2}=p^{\perp}\mathcal{H}$ are orthogonal closed subspaces of $\mathcal{H}$. We call the projection $p$ on $\mathcal{H}$ nontrivial if $p$ is not equal to $0$ and $\id_{\mathcal{H}}$.

Using the projections on $\mathcal{H}$, we can construct Rota-Baxter operators on the $C^{\ast}$-algebra $B(\mathcal{H})$.

\begin{eg}\label{eg2}
Let $p$ be a projection on $\mathcal{H}$. We define a linear operator $L_{p}$ on $B(\mathcal{H})$ by
$$L_{p}(a)=pa$$
for any $a\in B(\mathcal{H})$. It is easy to verify that $L_{p}$ is a Rota-Baxter operator of weight $-1$ on $B(\mathcal{H})$. Then $(B(\mathcal{H}),L_{p})$ is a Rota-Baxter $C^{\ast}$-algebra.
\end{eg}

\subsection{Rota-Baxter operators on $B(\mathcal{H})$}\label{3.1}

From here to the end of this subsection, without further mention, we will always assume that $\mathcal{H}$ is a Hilbert space, $p$ is a projection on $\mathcal{H}$ with $\mathcal{H}_{1}=p\mathcal{H}$ and $\mathcal{H}_{2}=p^{\perp}\mathcal{H}$, $\mathcal{A}$ is a $C^{\ast}$-subalgebra of $B(\mathcal{H})$ and $\mathcal{A}_{1}$ is a $C^{\ast}$-subalgebra of $\mathcal{A}$. Note that we have $\mathcal{H}=\mathcal{H}_1\oplus\mathcal{H}_2$.

\begin{Def}
Let $f$ be a linear operator on $\mathcal{H}$ and let $P$ be a Rota-Baxter operator of weight $\lambda\in\mathbb{C}$ on $\mathcal{A}$.  The operator $P$ is called a Rota-Baxter operator matching $f$ on $\mathcal{A}_{1}$ if
$$P(a).f(x)=f(P(a).x)+f(a.f(x))+\lambda f(a.x),\qquad\forall a\in \mathcal{A}_{1},\ \forall x\in\mathcal{H}.$$
\end{Def}

As a special case of the above definition,  Rota-Baxter operators matching projections on the Hilbert space $\mathcal{H}$ are the central notion of the article. We first give an equivalent condition for a Rota-Baxter operator of weight $-1$ which matching projections. Recall that for a Rota-Baxter operator $P$ of weight $\lambda$ on $\mathcal{A}$, we set $\widetilde{P}=-\lambda \id_{\mathcal{A}}-P$.

\begin{lem}\label{lem match}
Let $P$ be a Rota-Baxter operator of weight $-1$ on $\mathcal{A}$. Then $P$ is a Rota-Baxter operator matching $p$ on $\mathcal{A}_{1}$ if and only if for any $a\in \mathcal{A}_{1}$ and any $x\in\mathcal{H}$, we have $P(a).p(x)\in \mathcal{H}_1$ and $\widetilde{P}(a).p^{\perp}(x)\in \mathcal{H}_2$.	
\end{lem}

\proof
Assume that $P$ matches $p$ on $\mathcal{A}_{1}$. For any $x\in \mathcal{H}$ and any $a\in\mathcal{A}_{1}$, we have
$$P(a).p^2(x)=p(P(a).p(x))+p(a.p^2(x))-p(a.p(x)).$$
As $p^2=p$, we have
$$P(a).p(x)=p(P(a).p(x))\in \mathcal{H}_1.$$
Similarly, one can show that $\widetilde{P}(a).p^{\perp}(x)\in \mathcal{H}_2$.

Conversely, assume that for any $a\in\mathcal{A}_1$ and any $x\in\mathcal{H}$, we have $P(a).p(x)\in \mathcal{H}_1$ and $\widetilde{P}(a).p^{\perp}(x)\in \mathcal{H}_2$. We write $x=x_{1}+x_{2}$ with $x_{1}=p(x)\in \mathcal{H}_1$ and $x_2=p^{\perp}(x)\in \mathcal{H}_2$. Then we have
\begin{align*}
&p(P(a).x)+p(a.p(x))-p(a.x)\\
=&p(P(a).(x_{1}+x_{2}))+p(a.x_{1})-p(a.(x_{1}+x_{2}))\\
=&p(P(a).x_{1})+p(P(a).x_{2})-p(a.x_{2})\\
=&p(P(a).x_{1})-p(\widetilde{P}(a).x_{2}).
\end{align*}
Since $P(a).x_1\in \mathcal{H}_1$ and $\widetilde{P}(a).x_2\in \mathcal{H}_2$, we have
$$p(P(a).x)+p(a.p(x))-p(a.x)=P(a).x_{1}=P(a).p(x),$$
which shows that $P$ matches $p$ on $\mathcal{A}_{1}$.
\qed

As $p^{\perp}$ is also a projection, the above lemma has an immediate corollary.

\begin{cor}
Let $P$ be a Rota-Baxter operator of weight $-1$ on $\mathcal{A}$. Then $P$ matches $p$ on $\mathcal{A}_{1}$ if and only if $\widetilde{P}$ matches $p^{\perp}$ on $\mathcal{A}_{1}$.
\end{cor}

Now we write an element $x=x_1+x_2\in\mathcal{H}$ with $x_1=p(x)\in\mathcal{H}_1$ and $x_2=p^{\perp}(x)\in\mathcal{H}_2$ as a vector $x=\begin{bmatrix}
x_1\\
x_2
\end{bmatrix}$. As shown in [\citen{N.P},Theorem 2.10], for any $a\in B(\mathcal{H})$, since
\begin{align*}
a.x&=a(p(x))+a(p^{\perp}(x))\\
=&p(a(p(x)))+p^{\perp}(a(p(x)))+p(a(p^{\perp}(x)))+p^{\perp}(a(p^{\perp}(x)))\\
=&(pap)(x_1)+(pap^{\perp})(x_2)+(p^{\perp}ap)(x_1)+(p^{\perp}ap^{\perp})(x_2),
\end{align*}
we may write $a$ as a matrix
\begin{equation*}
a=\begin{bmatrix}
a_{11} & a_{12}\\
a_{21} & a_{22}
\end{bmatrix},
\end{equation*}
with
\begin{align*}
&a_{11}=pap:\mathcal{H}_1\rightarrow\mathcal{H}_1,\qquad a_{12}=pap^{\perp}:\mathcal{H}_2\rightarrow\mathcal{H}_1, \\
&a_{21}=p^{\perp}ap:\mathcal{H}_1\rightarrow\mathcal{H}_2,\quad\;  a_{22}=p^{\perp}ap^{\perp}:\mathcal{H}_2\rightarrow\mathcal{H}_2.
\end{align*}
It is easy to show that the operations in $B(\mathcal{H})$ coincides with the corresponding matrix operations. For example, we have
$$a^{\ast}=\begin{bmatrix}
a_{11}^{\ast} & a_{21}^{\ast}\\
a_{12}^{\ast} & a_{22}^{\ast}
\end{bmatrix},$$
which comes from
\begin{align*}
&(a^{\ast})_{11}=pa^{\ast}p=(pap)^{\ast}=a_{11}^{\ast},\\
&(a^{\ast})_{12}=pa^{\ast}p^{\perp}=(p^{\perp}ap)^{\ast}=a_{21}^{\ast},\\
&(a^{\ast})_{21}=p^{\perp}a^{\ast}p=(pap^{\perp})^{\ast}=a_{12}^{\ast},\\
&(a^{\ast})_{22}=p^{\perp}a^{\ast}p^{\perp}=(p^{\perp}ap^{\perp})^{\ast}=a_{22}^{\ast}.
\end{align*}
And if $b=\begin{bmatrix}
b_{11} & b_{12}\\
b_{21} & b_{22}
\end{bmatrix}\in B(\mathcal{H})$, then
$$ab=\begin{bmatrix}
a_{11}b_{11}+a_{12}b_{21} & a_{11}b_{12}+a_{12}b_{22}\\
a_{21}b_{11}+a_{22}b_{21} & a_{21}b_{12}+a_{22}b_{22}
\end{bmatrix}.$$
In fact, we have
\begin{align*}
	ab=&(pap+pap^{\perp}+p^{\perp}ap+p^{\perp}ap^{\perp})(pbp+pbp^{\perp}+p^{\perp}bp+p^{\perp}bp^{\perp})\\
	=&pappbp+pap^{\perp}p^{\perp}bp+pappbp^{\perp}+pap^{\perp}p^{\perp}bp^{\perp}+p^{\perp}appbp\\
	&+p^{\perp}ap^{\perp}p^{\perp}bp+p^{\perp}appbp^{\perp}+p^{\perp}ap^{\perp}p^{\perp}bp^{\perp}\\
	=&p(a_{11}b_{11}+a_{12}b_{21})p+p(a_{11}b_{12}+a_{12}b_{22})p^{\perp}+p^{\perp}(a_{21}b_{11}+a_{22}b_{21})p\\
	&+p^{\perp}(a_{21}b_{12}+a_{22}b_{22})p^{\perp}.
\end{align*}
Using the above matrix notation, we may restate Lemma \ref{lem match} as in the following proposition.

\begin{prop}\label{prop divide matrix}
Let $P$ be a Rota-Baxter operator of weight $-1$ on $\mathcal{A}$. Then $P$ matches $p$ on $\mathcal{A}_{1}$ if and only if for any $a\in\mathcal{A}_{1}$, we have $P(a)_{21}=0$ and $P(a)_{12}=a_{12}$.
\end{prop}

\proof
 Assume that $P$ matches $p$ on $\mathcal{A}_1$. For any $x_1\in\mathcal{H}_1$, we have
 $$P(a)_{21}(x_1)=(p^{\perp}P(a)p)(p(x))=p^{\perp}(P(a).p(x)).$$
By Lemma \ref{lem match}, $P(a).p(x)\in \mathcal{H}_1$. Hence $P(a)_{21}(x_1)=0$ and then $P(a)_{21}=0$. Similarly, for any $x_2\in\mathcal{H}_2$, we have
$$P(a)_{12}(x_2)=(pP(a)p^{\perp})(p^{\perp}(x))=p(P(a).p^{\perp}(x)).$$
Since $\widetilde{P}=\id-P$, we have
$$P(a)_{12}(x_2)=p(a.p^{\perp}(x))-p(\widetilde{P}(a).p^{\perp}(x))=a_{12}(x_2)-p(\widetilde{P}(a).p^{\perp}(x)).$$
By Lemma  \ref{lem match}, $\widetilde{P}(a).p^{\perp}(x)\in \mathcal{H}_2$. Therefore $P(a)_{12}(x_2)=a_{12}(x_2)$ and then $P(a)_{12}=a_{12}$.

Conversely, assume that for any $a\in \mathcal{A}_{1}$, we have $P(a)_{21}=0$ and $P(a)_{12}=a_{12}$. For any $x\in\mathcal{H}$, we have
\begin{align*}
P(a).p(x)=&pP(a)p(x_1)+p^{\perp}P(a)p(x_1)\\
=&pP(a)p(x_1)+P(a)_{21}(x_1)=pP(a)p(x_1)\in \mathcal{H}_1,
\end{align*}
and
\begin{align*}
\widetilde{P}(a).p^{\perp}(x)=&p\widetilde{P}(a)p^{\perp}(x_2)+p^{\perp}\widetilde{P}(a)p^{\perp}(x_2)\\
=&pap^{\perp}(x_2)-pP(a)p^{\perp}(x_2)+p^{\perp}\widetilde{P}(a)p^{\perp}(x_2)\\
=&(a_{12}-P(a)_{12})(x_2)+p^{\perp}\widetilde{P}(a)p^{\perp}(x_2)=p^{\perp}\widetilde{P}(a)p^{\perp}(x_2)\in \mathcal{H}_2.
\end{align*}
Hence by Lemma \ref{lem match}, the Rota-Baxter operator $P$ matches $p$ on $\mathcal{A}_1$.
\qed

Note that the direct sum $B(\mathcal{H}_1)\oplus B(\mathcal{H}_2)$ can be regarded as a $C^{\ast}$-subalgebra of $B(\mathcal{H})$ by the following embedding
$$\iota:B(\mathcal{H}_1)\oplus B(\mathcal{H}_2)\rightarrow B(\mathcal{H}); (a_1,a_2)\mapsto\begin{bmatrix}
a_1 & 0\\
0 & a_2
\end{bmatrix}.$$
The following is the main result of this note.

\begin{thm}\label{th}
Let $\mathcal{H}$ be a Hilbert space and let $p$ be a projection on $\mathcal{H}$. Set $\mathcal{H}_{1}=p\mathcal{H}$ and $\mathcal{H}_{2}=p^{\perp}\mathcal{H}$. Then there is a $1$-$1$ correspondence between the following two sets:
\begin{description}
\item[(1)] the set of Rota-Baxter operators of weight $-1$ on $B(\mathcal{H})$ matching $p$ on $B(\mathcal{H})$;
\item[(2)] the set of Rota-Baxter operators of weight $-1$ on $B(\mathcal{H}_{1})\oplus B(\mathcal{H}_{2})$.
\end{description}
\end{thm}

\proof
Denote by $\Gamma$ the set of Rota-Baxter operators of weight $-1$ on $B(\mathcal{H})$ matching $p$ on $B(\mathcal{H})$, and by $\Lambda$ the set of Rota-Baxter operators of weight $-1$ on $B(\mathcal{H}_{1})\oplus B(\mathcal{H}_{2})$.

Assume that $P\in \Gamma$. For any $a_{1}\in B(\mathcal{H}_{1})$ and $a_{2}\in B(\mathcal{H}_{2})$, we have $a=\begin{bmatrix}
a_1 & 0\\
0 & a_2
\end{bmatrix}\in B(\mathcal{H})$ and
$$P(a)=\begin{bmatrix}
P(a)_{11} & 0\\
0 & P(a)_{22}
\end{bmatrix}\in B(\mathcal{H}_{1})\oplus B(\mathcal{H}_{2})$$
by Proposition \ref{prop divide matrix}. Hence we get a linear operator $P'$ on $B(\mathcal{H}_{1})\oplus B(\mathcal{H}_{2})$ defined by $P'=P|_{B(\mathcal{H}_{1})\oplus B(\mathcal{H}_{2})}$. It is obvious  that $P'$  is a Rota-Baxter operator of weight $-1$ on $B(\mathcal{H}_{1})\oplus B(\mathcal{H}_{2})$. Hence we have a map $\Phi:\Gamma\rightarrow\Lambda$ defined by $\Phi(P)=P'$.

Conversely,  assume that $P'\in\Lambda$. For any $a=\begin{bmatrix}
a_{11} & a_{12}\\
a_{21} & a_{22}
\end{bmatrix}\in B(\mathcal{H})$, we have $\begin{bmatrix}
a_{11} & 0\\
0 & a_{22}
\end{bmatrix}\in B(\mathcal{H}_{1})\oplus B(\mathcal{H}_{2})$. We then define an operator $P$ on $B(\mathcal{H})$ by
$$P(a)=P'\left(\begin{bmatrix}
a_{11} & 0  \\
0 & a_{22}
\end{bmatrix}\right)+\begin{bmatrix}
0 &a_{12}\\
0 & 0
\end{bmatrix}.$$
It is obvious that $P$ is linear.

We show that $P$ is a Rota-Baxter operator of weight $-1$ on $B(\mathcal{H})$. Assume that $P'\left(\begin{bmatrix}
a_{11} & 0  \\
0 & a_{22}
\end{bmatrix}\right)=
\begin{bmatrix}
a'_{11} & 0  \\
0 & a'_{22}
\end{bmatrix}$, then we have
$$P(a)=
\begin{bmatrix}
a'_{11} & 0  \\
0 & a'_{22}
\end{bmatrix}+\begin{bmatrix}
0 &a_{12}\\
0 & 0
\end{bmatrix}.$$
For another element $b=\begin{bmatrix}
b_{11} &b_{12}\\
b_{21} &b_{22}
\end{bmatrix}\in B(\mathcal{H})$, we set $P'\left(\begin{bmatrix}
b_{11} & 0  \\
0 & b_{22}
\end{bmatrix}\right)=
\begin{bmatrix}
b'_{11} & 0  \\
0 & b'_{22}
\end{bmatrix}$. Hence we have
$$P(a)P(b)=\begin{bmatrix}
		a'_{11} &a_{12}\\
		0 & a'_{22}
	\end{bmatrix}\begin{bmatrix}
		b'_{11} &b_{12}\\
		0 & b'_{22}
	\end{bmatrix}=\begin{bmatrix}
		a'_{11}b'_{11} & a'_{11}b_{12}+a_{12}b'_{22}  \\
		0 & a'_{22}b'_{22}
	\end{bmatrix}.$$
Since
$$aP(b)=\begin{bmatrix}
a_{11} & a_{12}  \\
a_{21} & a_{22}
\end{bmatrix}
\begin{bmatrix}
b'_{11} & b_{12}  \\
0 & b'_{22}
\end{bmatrix}=\begin{bmatrix}
a_{11}b'_{11} & a_{11}b_{12}+a_{12}b'_{22}  \\
a_{21}b'_{11} & a_{21}b_{12}+a_{22}b'_{22}
\end{bmatrix},$$
we get
$$P(aP(b))=P'\left(\begin{bmatrix}
a_{11}b'_{11} & 0  \\
0 & a_{21}b_{12}+a_{22}b'_{22}
\end{bmatrix}\right)+\begin{bmatrix}
0 & a_{11}b_{12}+a_{12}b'_{22}  \\
0 & 0
\end{bmatrix}.$$
Similarly, we have
$$P(P(a)b)=P'\left(\begin{bmatrix}
a'_{11}b_{11}+a_{12}b_{21} & 0  \\
0 & a'_{22}b_{22}
\end{bmatrix}\right)+\begin{bmatrix}
0 & a'_{11}b_{12}+a_{12}b_{22}  \\
0 & 0
\end{bmatrix}$$
and
$$P(ab)=P'\left(\begin{bmatrix}
a_{11}b_{11}+a_{12}b_{21} & 0  \\
0 & a_{21}b_{12}+a_{22}b_{22}
\end{bmatrix}\right)+\begin{bmatrix}
0 & a_{11}b_{12}+a_{12}b_{22}  \\
0 & 0
\end{bmatrix}.$$
Therefore we find
\begin{align*}
&P(aP(b))+P(P(a)b)-P(ab)\\
=&P'\left(\begin{bmatrix}
a_{11}b'_{11}+a'_{11}b_{11}-a_{11}b_{11} & 0  \\
0 & a_{22}b'_{22}+a'_{22}b_{22}-a_{22}b_{22}
\end{bmatrix}\right)+\begin{bmatrix}
0 & a_{12}b'_{22}+a'_{11}b_{12}  \\
0 & 0
\end{bmatrix}.
\end{align*}
On the other hand, since $P'\in\Lambda$, we have
\begin{align*}
&\begin{bmatrix}
   a'_{11} & 0\\
   0 & a'_{22}
   \end{bmatrix}\begin{bmatrix}
   b'_{11} & 0\\
   0 & b'_{22}
   \end{bmatrix}\\
   =&P'\left(\begin{bmatrix}
   a_{11} & 0\\
   0 & a_{22}
   \end{bmatrix}\begin{bmatrix}
   b'_{11} & 0\\
   0 & b'_{22}
   \end{bmatrix}+\begin{bmatrix}
   a'_{11} & 0\\
   0 & a'_{22}
   \end{bmatrix}\begin{bmatrix}
   b_{11} & 0\\
   0 & b_{22}
   \end{bmatrix}-\begin{bmatrix}
   a_{11}b_{11} & 0\\
   0 & a_{22}b_{22}
   \end{bmatrix}\right).
\end{align*}
Then we get
$$P(aP(b))+P(P(a)b)-P(ab)=\begin{bmatrix}
   a'_{11}b'_{11} & 0\\
   0 & a'_{22}b'_{22}
   \end{bmatrix}+\begin{bmatrix}
0 & a_{12}b'_{22}+a'_{11}b_{12}  \\
0 & 0
\end{bmatrix},$$
which implies that
$$P(a)P(b)=P(aP(b))+P(P(a)b)-P(ab).$$
So $P$ is a Rota-Baxter operator of weight $-1$ on $B(\mathcal{H})$. And from Proposition \ref{prop divide matrix}, we see that $P$ matches $p$ on $B(\mathcal{H})$. Then we have a map $\Psi:\Lambda\rightarrow\Gamma$ defined by $\Psi(P')=P$.

It is easy to verify that $\Psi\circ\Phi$ and $\Phi\circ\Psi$ are identity maps. And so there is a $1$-$1$ correspondence between $\Lambda$ and $\Gamma$.
\qed

\begin{rem}
We can't extend the above theorem to the $C^{\ast}$-subalgebra $\mathcal{A}$ of $B(\mathcal{H})$. For a Rota-Baxter operator $P'$ of weight $-1$ on $p\mathcal{A}p\oplus p^{\perp}\mathcal{A}p^{\perp}$, as $p$ and $p^{\perp}$ may not be in $\mathcal{A}$, we can't construct a Rota-Baxter operator  $P$ of weight $-1$ matching $p$ on $\mathcal{A}$ as above.
\end{rem}

From Theorem \ref{th}, we can construct Rota-Baxter operators of weight $-1$ on $B(\mathcal{H})$  matching $p$ on $B(\mathcal{H})$ from the Rota-Baxter operators of weight $-1$ on $B(\mathcal{H}_1)\oplus B(\mathcal{H}_2)$. The following lemma gives a method to construct  Rota-Baxter operators on $B(\mathcal{H}_1)\oplus B(\mathcal{H}_2)$.

\begin{lem}\label{lem direct sum}
Let $\mathcal{A}_{1}$ and $\mathcal{A}_{2}$ be two $C^{\ast}$-algebras. If $P_{1}$ and $P_{2}$ are Rota-Baxter operators of weight $\lambda\in\mathbb{C}$ on $\mathcal{A}_{1}$ and $\mathcal{A}_{2}$, respectively, then the operator
$$P=P_{1}\oplus P_{2}:\mathcal{A}_{1} \oplus\mathcal{A}_{2}\rightarrow \mathcal{A}_{1} \oplus\mathcal{A}_{2}; (a_1,a_2)\mapsto (P_1(a_1),P_2(a_2))$$
is a Rota-Baxter operator  of weight $\lambda$ on $\mathcal{A}_{1} \oplus\mathcal{A}_{2}$.
\end{lem}

\proof Is is trivial to check.\qed

Then using Lemma \ref{lem direct sum} and Theorem \ref{th}, we get

\begin{cor}\label{3.8}
Let $P_{1}$ be a Rota-Baxter operator of weight $-1$ on $B(\mathcal{H}_1)$ and let $P_{2}$ be a Rota-Baxter operator of weight $-1$ on $B(\mathcal{H}_2)$. Then the operator
$$P=\begin{bmatrix}
  		P_{1}  &   \id\\
  		0     &   P_{2}\\  		
  	\end{bmatrix}:B(\mathcal{H})\rightarrow B(\mathcal{H});\begin{bmatrix}
  		a_{11}  &   a_{12}\\
  		a_{21}     &   a_{22}\\  		
  	\end{bmatrix}\mapsto \begin{bmatrix}
  		P_1(a_{11})  &   a_{12}\\
  		0     &   P_2(a_{22})\\  		
  	\end{bmatrix}	$$
is a Rota-Baxter operator of weight $-1$ matching $p$ on $B(\mathcal{H})$.
\end{cor}

Recall that a Rota-Baxter operator $P$ of weight $\lambda$ on $\mathcal{A}$  is called symmetric on $\mathcal{A}_{1}$ provided that
$$P(a^{\ast})=P(a)^{\ast}, \qquad \forall a\in \mathcal{A}_{1}.$$

\begin{prop}\label{sym}
If there is a Rota-Baxter operator $P$ of weight $-1$ on $\mathcal{A}$ which is symmetric and matches $p$ on $\mathcal{A}_{1}$, then we have $\mathcal{A}_{1}\subseteq B(\mathcal{H}_1)\oplus B(\mathcal{H}_2)$.
\end{prop}

\proof
For any $a=\begin{bmatrix}
		a_{11}   &   a_{12}\\
  		a_{21}   &   a_{22}
	\end{bmatrix}\in \mathcal{A}_{1}$, from Proposition \ref{prop divide matrix}, we have
$$P(a)=\begin{bmatrix}
		P(a)_{11}   &   a_{12}\\
		0      &   P(a)_{22}
	\end{bmatrix}.$$
On the other hand, since $P$ is symmetric on $\mathcal{A}_{1}$, we have
$$P(a)=P(a^{\ast})^{\ast}
   =\begin{bmatrix}
   		(P(a^{\ast})_{11})^{\ast}   &   0\\
   	         a_{21}                         &   (P(a^{\ast})_{22})^{\ast}
            \end{bmatrix}.$$
Then we get $a_{12}=0$ and $a_{21}=0$. So $\mathcal{A}_{1}\subseteq B(\mathcal{H}_1)\oplus B(\mathcal{H}_2)$.
\qed

To state the last result of this subsection, we recall the definition of $\ast$-representation  from \cite{Sakai1}. A $\ast$-representation $\pi$ of the $C^{\ast}$-algebra $\mathcal{A}$ is a $\ast$-homomorphism $\pi: \mathcal{A}\rightarrow B(\mathcal{H})$ with $\mathcal{H}$ a Hilbert space. We denote this $\ast$-representation of $\mathcal{A}$ by $\{\pi,\mathcal{H}\}$. We say $\{\pi,\mathcal{H}\}$ is topologically irreducible if $\pi(\mathcal{A})$ has no proper invariant subspaces. It is called algebraically irreducible if it has no proper invariant manifolds(subspaces of $\mathcal{H}$ that are not necessarily closed). These two notions coincide for $C^{\ast}$-algebras from \cite{Davidson}. Hence we call $\{\pi,\mathcal{H}\}$ irreducible when either of the two conditions holds.

\begin{cor}\label{3.11}
Let $\{\pi,\mathcal{H}\}$ be a $\ast$-representation of a $C^{\ast}$-algebra $\mathcal{A}$. Then the representation $\{\pi,\mathcal{H}\}$ is irreducible if and only if  there is no Rota-Baxter operator of weight $-1$ on $B(\mathcal{H})$ which is symmetric and matches $p$ on $\pi(\mathcal{A})$ for some nontrivial projection $p$ on $\mathcal{H}$.
\end{cor}

\proof
If there is such a Rota-Baxter operator on $\pi(\mathcal{A})$, then from Proposition \ref{sym}, we have $\pi(\mathcal{A})\subseteq B(p\mathcal{H})\oplus B(p^{\perp}\mathcal{H})$. Hence $\pi(\mathcal{A})$ has a proper invariant closed subspace $p\mathcal{H}$  of $\mathcal{H}$. Therefore $\{\pi,\mathcal{H}\}$ is not irreducible.

Conversely, if the representation $\{\pi,\mathcal{H}\}$ is not irreducible, then $\pi(\mathcal{A})$ has a proper invariant closed subspace $\mathcal{H}_{1}$. There exists a projection $p$ on $\mathcal{H}$ such that $\mathcal{H}_{1}=p\mathcal{H}$. Therefore $\pi(\mathcal{A}).p\mathcal{H}\subseteq p\mathcal{H}$. Define a linear operator $P$ on $B(\mathcal{H})$ by
$$P(a)=pap+p^{\perp}ap^{\perp}+pap^{\perp}$$
for any $a\in B(\mathcal{H})$. Then by the proof of Theorem \ref{th}, $P$ is a Rota-Baxter operator of weight $-1$ matching $p$ on $B(\mathcal{H})$. For any $b\in \pi(\mathcal{A})$ , we have $b=P(b)$ and $b^{\ast}=P(b^{\ast})$ from Proposition \ref{prop divide matrix}. Hence we find
$$P(b)=b=(b^{\ast})^{\ast}=P(b^{\ast})^{\ast},$$
which means that $P$ is symmetric on $\pi(\mathcal{A})$.
\qed

At the end of this subsection, we see an example.

\begin{eg}
	Let $f:\mathcal{A}\to\mathbb{C}$ be a linear function, where $\mathcal{A}$ is a $C^{\ast}$-algebra. The function $f$ is called a positive linear functional if $f(a)\ge0$ for any $a\ge0$ in $\mathcal{A}$ (the notion of $a\ge 0$ can be found in Section $\ref{4}$). The norm of $f$ is defined by
$$\|f\|=\underset{a\in\mathcal{A},a\neq 0}{\sup}\frac{|f(a)|}{\|a\|}.$$
As in \cite{Davidson}, the function $f$ is called a state on $\mathcal{A}$ if $f$ is a positive linear functional and $\|f\|=1$. The state $f$ is called a pure state if there is not a $\lambda\in (0,1)$ and two states $f_{1}$ and $f_{2}$ on $\mathcal{A}$, such that $f=\lambda f_{1}+(1-\lambda)f_{2}$.

	From \cite[Theorem \uppercase\expandafter{\romannumeral1}.9.6]{Davidson}, for any state $f$ on $\mathcal{A}$, we can construct a Hilbert space $\mathcal{H}_{f}$ and a $\ast$-representation $\{\pi_{f},\mathcal{H}_{f}\}$ of $\mathcal{A}$. This method of building representations from states is called GNS-constructions. From \cite[Theorem \uppercase\expandafter{\romannumeral1}.9.8]{Davidson}, we know that the representation $\{\pi_{f},\mathcal{H}_{f}\}$ is an irreducible representation of $\mathcal{A}$ if and only if $f$ is a pure state. Hence if $f$ is not a pure state, then from Corollary \ref{3.11}, there is at least one Rota-Baxter operator of weight $-1$ which is symmetric and matches $p$ on $\pi_{f}(\mathcal{A})$ for some nontrivial projection $p$ on $\mathcal{H}_f$.
\end{eg}

\subsection{Representations of Rota-Baxter $C^{\ast}$-algebras}\label{3.2}

We introduce the notion of representations of Rota-Baxter $C^{\ast}$-algebras.

\begin{Def}
Let $\mathcal{A}$ be a $C^{\ast}$-algebra and let $P$ be a Rota-Baxter operator of weight $\lambda$ on $\mathcal{A}$. Let $\mathcal{H}$ be a Hilbert space and let $P'$ be a Rota-Baxter operator of weight $\lambda$ on $B(\mathcal{H})$.  Let $\pi:(\mathcal{A},P)\rightarrow (B(\mathcal{H}),P')$ be a Rota-Baxter $C^\ast$-algebra homomorphism and let $f$ be a linear operator on $\mathcal{H}$. Then $\{\pi,\mathcal{H},f\}$ is a Rota-Baxter $\ast$-representation of $(\mathcal{A},P)$ into $(B(\mathcal{H}),P')$ if $P'$ matches $f$ on $\pi(\mathcal{A})$.
\end{Def}

Hence $\{\pi,\mathcal{H},f\}$ is a Rota-Baxter $\ast$-representation of $(\mathcal{A},P)$ into $(B(\mathcal{H}),P')$ if for any $h\in\mathcal{H}$ and any $a\in\mathcal{A}$ we have
 $$\pi (P(a))(f(h))=f(\pi (P(a))(h))+f(\pi (a)(f(h)))+\lambda f(\pi (a)(h)).$$ 	

We first prepare two lemmas.

\begin{lem}\label{rep1}
Let $\mathcal{H}$ be a Hilbert space and let $p$ be a projection on $\mathcal{H}$. Let $P'$ be a Rota-Baxter operator of weight $-1$ on $B(p\mathcal{H})\oplus B(p^{\perp}\mathcal{H})$ and let $P$ be the Rota-Baxter operator  of weight $-1$ on $B(\mathcal{H})$ matching $p$ that corresponding to $P'$ (See Theorem \ref{th}.). Then the embedding $\iota:B(p\mathcal{H})\oplus B(p^{\perp}\mathcal{H})\rightarrow B(\mathcal{H})$ is a Rota-Baxter $C^\ast$-algebra homomorphism.
\end{lem}

\proof
It needs to show that $\iota\circ P'=P\circ \iota$, which is obvious.
\qed

\begin{lem}\label{rep2}
Let  $\pi_{1}:(\mathcal{A}_{1},P_{1})\to(\mathcal{A}_{2},P_{2})$ and $\pi_{2}:(\mathcal{A}'_{1},P'_{1})\to(\mathcal{A}'_{2},P'_{2})$ be two Rota-Baxter $C^\ast$-algebra homomorphisms. Then the map
$$\pi=\pi_1\oplus\pi_2:\mathcal{A}_{1}\oplus \mathcal{A}'_{1}\rightarrow \mathcal{A}_{2}\oplus \mathcal{A}'_{2};(a_1,b_1)\mapsto (\pi_1(a_1),\pi_2(b_1))$$
is a Rota-Baxter $C^\ast$-algebra homomorphism.
\end{lem}

\proof
It is easy to verify that $\pi$ is a Rota-Baxter $C^\ast$-algebra homomorphism.
\qed

Then we state the main theorem of this subsection.

\begin{thm}
Let $\mathcal{H}_{1}$ and $\mathcal{H}_{2}$ be two Hilbert spaces. For $i=1,2$, let $P'_i$ be a Rota-Baxter operator of weight $-1$ on  $B(\mathcal{H}_{i})$, $\mathcal{A}_{i}$ be a $C^{\ast}$-subalgebra of  $B(\mathcal{H}_{i})$, and let $P_{i}$ be a Rota-Baxter operator of weight $-1$ on  $\mathcal{A}_{i}$. Set $\mathcal{H}=\mathcal{H}_{1}\oplus \mathcal{H}_{2}$.  Let $P'$ be the unique Rota-Baxter operator of weight $-1$ on $B(\mathcal{H})$ corresponding to $P'_{1}\oplus P'_{2}$ on $B(\mathcal{H}_{1})\oplus B(\mathcal{H}_{2})$ (See Theorem \ref{th} and Lemma \ref{lem direct sum}.).  Let $p$ be the projection from $\mathcal{H}$ onto $\mathcal{H}_{1}$. Then the following two statements are equivalent:
\begin{description}
    	\item[(a)] There are Rota-Baxter $\ast$-representations $\{\pi_{i},\mathcal{H}_i,\id_{\mathcal{H}_i}\}$ from $(\mathcal{A}_{i},P_{i})$ into $(B(\mathcal{H}_{i}),P'_{i})$ with $i=1,2$;
    	\item[(b)] There is a Rota-Baxter $\ast$-representation $\{\pi,\mathcal{H},p\}$ from $(\mathcal{A}_{1}\oplus \mathcal{A}_{2},P_{1}\oplus P_{2})$ into $(B(\mathcal{H}), P')$ such that $\pi(\mathcal{A}_{i})\subseteq B(\mathcal{H}_{i})$ for $i=1,2$.
\end{description}
\end{thm}

\proof
Set $\mathcal{A}=\mathcal{A}_{1}\oplus \mathcal{A}_{2}$ and $P=P_{1}\oplus P_{2}$.

(a) $\Rightarrow$ (b) Assume that (a) holds. Using Lemmas \ref{rep1} and \ref{rep2}, we have a Rota-Baxter $C^\ast$-algebra homomorphism
  $$\pi:(\mathcal{A}_{1}\oplus \mathcal{A}_{2},P_{1}\oplus P_{2})\overset{\pi_1\oplus\pi_2}{\longrightarrow} (B(\mathcal{H}_1)\oplus B(\mathcal{H}_2),P'_1\oplus P'_2)\overset{\iota}{\longrightarrow} (B(\mathcal{H}), P').$$
It is obvious that $\pi(\mathcal{A}_{i})\subseteq B(\mathcal{H}_{i})$ for $i=1,2$.

We check that $P'$ matches $p$ on $\pi(\mathcal{A})$. For any $a=(a_1,a_2)\in\mathcal{A}$ with $a_i\in\mathcal{A}_i$ and any $h=(h_1,h_2)\in\mathcal{H}$ with $h_i\in\mathcal{H}_i$, we have
\begin{align*}
&p(\pi (P(a))(h))+p(\pi (a)(p(h)))-p(\pi (a)(h))\\
=&p(\pi_1(P_1(a_1))(h_1),\pi_2(P_2(a_2))(h_2))+p(\pi_1(a_1)(h_1),0)-p(\pi_1(a_1)(h_1),\pi_2(a_2)(h_2))\\
=&\pi_1(P_1(a_1))(h_1),
\end{align*}
and
$$\pi (P(a))(p(h))=(\pi_1(P_1(a_1)),\pi_2(P_2(a_2)))(h_1,0)=\pi_1(P_1(a_1))(h_1).$$
Hence $P'$ matches $p$ on $\pi(\mathcal{A})$.

(a) $\Leftarrow$ (b) Assume that (b) holds. As $\pi$ is a Rota-Baxter $C^\ast$-algebra homomorphism, for any $a_{1}\in \mathcal{A}_{1}$ and $a_{2}\in \mathcal{A}_{2}$, we have
  $$\pi \circ P(a_{1},a_{2})=P'\circ\pi(a_{1},a_{2}).$$
  Taking $a_{2}=0$, we get $\pi\circ P_1(a_{1})=P'\circ\pi(a_{1})$. Since $P_{1}(a_{1})\in \mathcal{A}_{1}$ and $\pi(a_{1})\in B(\mathcal{H}_{1})$, we have
  $$\pi\circ P_{1}(a_1)=P'_{1}\circ \pi(a_1).$$
  Hence we get a Rota-Baxter $C^\ast$-algebra homomorphism
  $$\pi_{1}=\pi|_{\mathcal{A}_{1}}:(\mathcal{A}_{1},P_{1})\rightarrow (B(\mathcal{H}_{1}),P'_{1}).$$
Similarly, we have a Rota-Baxter $C^\ast$-algebra homomorphism
$$\pi_{2}=\pi|_{\mathcal{A}_{2}}: (\mathcal{A}_{2},P_{2})\rightarrow (B(\mathcal{H}_{2}),P'_{2}).$$
Finally, since $P'$ matches $p$ on $\pi(\mathcal{A})$, one can show that $P'_i$ matches $\id_{\mathcal{H}_i}$ on $\pi_i(\mathcal{A}_i)$ for $i=1,2$.
\qed



\section{Quasidiagonal operators and Rota-Baxter operators}\label{4}

In this section, all Hilbert spaces and $C^{\ast}$-algebras are assumed to be separable.
 Here a topology space $X$ is separable if it contains a countable dense subset.


We recall some basic definitions of quasidiagonal linear operators from \cite{N.P}. Let $\mathcal{A}$ be a $C^{\ast}$-algebra and let $a$ be a self-adjoint element of $\mathcal{A}$. If $\mathcal{A}$ is unital, the spectrum of $a$ is defined by
$$\Spec(a)=\{\lambda\in\mathbb{C}|a-\lambda 1_{\mathcal{A}} \text{\;is not invertible}\},$$
where $1_{\mathcal{A}}$ is the identity element of $\mathcal{A}$. If $\mathcal{A}$ is nonunital, then from [\citen{Sakai1}, Proposition 1.1.7], there is a unital $C^{\ast}$-algebra $\mathcal{A}_{1}$ such that $\mathcal{A}_{1}\simeq\mathcal{A}\bigoplus\mathbb{C}$. Then for any $a\in\mathcal{A}$, the spectrum of $a$ is the spectrum of $a$ as an element of $\mathcal{A}_{1}$.
The element $a$ is called positive if $\Spec(a)$ is contained in the non-negative reals. For any $a,b\in\mathcal{A}$, we say that $a$ is not bigger than $b$ if $b-a$ is positive, and we denote this by $a\le b$. We also denote $ab-ba$ by the Lie bracket $[a,b]$.

Let $\mathcal{H}$ be a Hilbert space. A projection $p$ on $\mathcal{H}$ is of finite rank if $p\mathcal{H}$ is of finite dimension. For any sequence ${a_{n}}\in B(\mathcal{H})$ and $a\in B(\mathcal{H})$, we say that ${a_{n}}\to a$ as $n\to\infty$ in the strong operator topology if $\|a_{n}.h-a.h\|\to0$ as $n\to\infty$ for any $h\in\mathcal{H}$.

We first give the definition of block diagonal linear operators.

\begin{Def}
	A bounded linear operator $d$ on a Hilbert space $\mathcal{H}$ is called block diagonal if there exists an increasing sequence of finite rank projections $p_{1}\le p_{2}\le p_{3}\le\cdots$ on $\mathcal{H}$, such that $\|[d,p_{n}]\|=0$ and $p_{n}\to \id_{\mathcal{H}}$ as $n\to\infty$ in the strong operator topology.
\end{Def}

We can construct an equivalent condition of block diagonal operators with the help of Rota-Baxter operators. Let $\mathcal{A}$ be a $C^{\ast}$-algebra. For any $a\in\mathcal{A}$, there is a $C^{\ast}$-subalgebra of $\mathcal{A}$ generated by $a$, which is denoted by $C^{\ast}(a)$.  In fact, $C^{\ast}(a)$ is the intersection of all $C^{\ast}$-subalgebras of $\mathcal{A}$ containing $a$, and $C^{\ast}(a)$ is the closure of the algebra generated by $a$ and $a^{\ast}$.

\begin{prop}\label{Prop:BlockDiagonal}
	Let $\mathcal{H}$ be a Hilbert space and $d\in B(\mathcal{H})$. Then the following statements are equivalent:
	\begin{description}
		\item[(a)] The operator $d$ is block diagonal;
		\item[(b)] There exists an increasing sequence of finite rank projections $p_{1}\le p_{2}\le p_{3}\le\cdots$ on $\mathcal{H}$ and a sequence of Rota-Baxter operators $\{P_n\}$ of weight $-1$ on $B(\mathcal{H})$, such that $p_{n}\to \id_{\mathcal{H}}$ as $n\to\infty$ in the strong operator topology, and $P_n$ are symmetric on $C^{\ast}(d)$ and matches $p_n$ on $B(\mathcal{H})$ for any $n\ge 1$.
	\end{description}
\end{prop}

\proof
For any $x\in \mathcal{H},n\in \mathbb{N}$ and $b\in B(\mathcal{H})$, we have $p_{n}^{\perp}bp_{n}(x)\in p_{n}^{\perp}\mathcal{H}$ and $p_{n}bp_{n}^{\perp}(x)\in p_{n}\mathcal{H}$. Therefore we find
$$
\|[b,p_{n}]\|=\|p_{n}^{\perp}bp_{n}-p_{n}bp_{n}^{\perp}\|\ge\|p_{n}^{\perp}bp_{n}\|.$$
And similarly we have $$\|[b,p_{n}]\|\ge\|p_{n}bp_{n}^{\perp}\|.$$
On the other hand, it holds
$$\|[b,p_{n}]\|=\|p_{n}^{\perp}bp_{n}-p_{n}bp_{n}^{\perp}\|\le\|p_{n}^{\perp}bp_{n}\|+\|p_{n}bp_{n}^{\perp}\|.$$
Therefore $\|[b,p_{n}]\|=0$ if and only if $\|p_{n}^{\perp}bp_{n}\|=0$ and $\|p_{n}bp_{n}^{\perp}\|= 0$.

(a) $\Rightarrow$ (b)
 Assume that $d$ is block diagonal. Then there exists an increasing sequence of finite rank projections $p_{1}\le p_{2}\le p_{3}\le\cdots$ on $\mathcal{H}$, such that $\|[d,p_{n}]\|=0$ for any $n\in\mathbb{N}$.

 Let $A$ be the set of elements which are finite products of $d$ and $d^{\ast}$ in $C^{\ast}(d)$. For any $b\in A$, if
$b=b_{1}b_{2}\cdots b_{n}$ with $b_{k}=d$ or $d^{\ast}$ for $1\le k \le n$, we set $|b|=n$. We first prove that for any $b\in A$, $\|[b,p_{n}]\|= 0$ by induction on $|b|$. The case of $|b|=1$ is trivial. For $|b|\ge2$, we can assume that there is a $b_{1}\in A$ such that $|b_{1}|=|b|-1$, and $b=b_{1}d$ without loss of generality. We have
 $$
 \begin{aligned}
 	\|p_{n}^{\perp}bp_{n}\|
 	&=\|p_{n}^{\perp}b_{1}p_{n}d+p_{n}^{\perp}b_{1}[d,p_{n}]\|\\
 	&\le\|p_{n}^{\perp}b_{1}p_{n}d\|+\|p_{n}^{\perp}b_{1}[d,p_{n}]\|\\
 	&\le\|p_{n}^{\perp}b_{1}p_n\|\|d\|+\|p_{n}^{\perp}\|\|b_{1}\|\|[d,p_{n}]\|.
 \end{aligned}
 $$
 Hence using the inductive assumption $\|p_n^{\perp}b_{1}p_n\|=0$ and the fact $\|[d,p_n]\|=0$, we find $\|p_{n}^{\perp}bp_{n}\|=0$. Similarly we get $\|p_{n}bp_{n}^{\perp}\|= 0$. Therefore we have $\|[b,p_{n}]\|= 0$.

 Using the above result, it is easy to know that for any $b\in C^{\ast}(d)$,  we have $\|[b,p_{n}]\|= 0$. Define a linear operator $P_{n}$ on $B(\mathcal{H})$ such that for any $a\in B(\mathcal{H})$ we have
$$P_n(a)=p_{n}ap_{n}+p^{\perp}_{n}ap^{\perp}_{n}+p_{n}ap^{\perp}_{n}.$$
Then $P_{n}$ is a Rota-Baxter operator of weight $-1$ matching $p_{n}$ on $B(\mathcal{H})$ by Corollary \ref{3.8}. Finally, for any $b\in C^{\ast}(d)$, as
$$
\begin{aligned}
	&\|P_n(b)-P_n(b^{\ast})^{\ast}\|\\
	&=\|p_{n}bp_{n}+p^{\perp}_{n}bp^{\perp}_{n}+p_{n}bp^{\perp}_{n}-p_{n}bp_{n}-p^{\perp}_{n}bp^{\perp}_{n}-p^{\perp}_{n}bp_{n}\|\\
	&=\|p_{n}bp^{\perp}_{n}-p^{\perp}_{n}bp_{n}\|\\
	&\le\|p_{n}bp^{\perp}_{n}\|+\|p^{\perp}_{n}bp_{n}\|=0,
\end{aligned}
$$	
we have $P_n(b)=P_n(b^{\ast})^{\ast}$.

(b) $\Rightarrow$ (a) Assume that (b) holds. Then using Proposition \ref{sym}, we know that $C^{\ast}(d)\subseteq B(p_{n}\mathcal{H})\bigoplus B(p_{n}^{\perp}\mathcal{H})$ for any $n\in\mathbb{N}$. Hence $\|p_{n}dp_{n}^{\perp}\|=\|p_{n}^{\perp}dp_{n}\|=0.$ Therefore  $\|[d,p_n]\|=0$, and (a) holds.
\qed

Now we give the definition of  quasidiagonal linear operators.

\begin{Def}
A bounded linear operator $d$ on a Hilbert space $\mathcal{H}$ is called quasidiagonal if there exists an increasing sequence of finite rank projections $p_{1}\le p_{2}\le p_{3}\le\cdots$ on $\mathcal{H}$, such that $\|[d,p_{n}]\|\to0$ as $n\to\infty$ and $p_{n}\to \id_{\mathcal{H}}$ as $n\to\infty$ in the strong operator topology.
\end{Def}

As above, we have an equivalent condition of quasidiagonal operators in terms of Rota-Baxter operators.

\begin{prop}
	Let $\mathcal{H}$ be a Hilbert space and $d\in B(\mathcal{H})$. Then the following statements are equivalent:
	\begin{description}
		\item[(a)] The operator $d$ is quasidiagonal;
		\item[(b)] There exists an increasing sequence of finite rank projections $p_{1}\le p_{2}\le p_{3}\le\cdots$ on $\mathcal{H}$ and a sequence of Rota-Baxter operators $\{P_n\}$ of weight $-1$ on $B(\mathcal{H})$, such that $P_{n}$ matches $p_{n}$ on $B(\mathcal{H})$ for any $n\geq 1$, $p_{n}\to \id_{\mathcal{H}}$ and $\|P_{n}(b)-P_{n}(b^{\ast})^{\ast}\|\to 0$ for any $b\in C^{\ast}(d)$ as $n\to\infty$.			
	\end{description}
\end{prop}

\proof
For any $x\in \mathcal{H}$ and $b\in B(\mathcal{H})$, similarly as in the proof of Proposition \ref{Prop:BlockDiagonal}, we have
$$\|[b,p_{n}]\|\ge\|p_{n}^{\perp}bp_{n}\|,\quad \|[b,p_{n}]\|\ge\|p_{n}bp_{n}^{\perp}\|,$$
and
$$\|[b,p_{n}]\|\le\|p_{n}^{\perp}bp_{n}\|+\|p_{n}bp_{n}^{\perp}\|.$$
Therefore $\|[b,p_{n}]\|\to 0$ as $n\to\infty$ if and only if $\|p_{n}^{\perp}bp_{n}\|\to0$ and $\|p_{n}bp_{n}^{\perp}\|\to 0$ as $n\to\infty$.

 (a) $\Rightarrow$ (b) Assume that $d$ is quasidiagonal. Then there exists an increasing sequence of finite rank projections $p_{1}\le p_{2}\le p_{3}\le\cdots$ on $\mathcal{H}$, such that $\|[d,p_{n}]\|\to0$ and $p_{n}\to \id_{\mathcal{H}}$ as $n\to\infty$.

Let $A$ be the set of elements which are finite products of $d$ and $d^{\ast}$ in $C^{\ast}(d)$.
We first prove that for any $b\in A$, $\|[b,p_{n}]\|\to 0$ as $n\to\infty$ by induction on $|b|$. The case of $|b|=1$ is trivial. For $|b|\ge2$, we can assume that there is a $b_{1}\in A$ such that $|b_{1}|=|b|-1$, and $b=b_{1}d$ without loss of generality. We have
$$
\begin{aligned}
	\|p_{n}^{\perp}bp_{n}\|
		&\le\|p_{n}^{\perp}b_{1}p_n\|\|d\|+\|p_{n}^{\perp}\|\|b_{1}\|\|[d,p_{n}]\|.
\end{aligned}
$$
Hence using the inductive assumption, we find $\|p_{n}^{\perp}bp_{n}\|\to 0$ as $n\to\infty$.
Similarly we get $\|p_{n}bp_{n}^{\perp}\|\to 0$ as $n\to\infty$. Therefore we have $\|[b,p_{n}]\|\to 0$ as $n\to \infty$.

Using the above result, it is easy to know that for any $b\in C^{\ast}(d)$,  we have $\|[b,p_{n}]\|\to 0$ as $n\to\infty$.

Define a linear operator $P_{n}$ on $B(\mathcal{H})$ such that for any $a\in B(\mathcal{H})$ we have
$$P_n(a)=p_{n}ap_{n}+p^{\perp}_{n}ap^{\perp}_{n}+p_{n}ap^{\perp}_{n}.$$
Then $P_{n}$ is a Rota-Baxter operator of weight $-1$ matching $p_{n}$ on $B(\mathcal{H})$ by Corollary \ref{3.8}. Finally, for any $b\in C^{\ast}(d)$, as
$$\|P_n(b)-P_n(b^{\ast})^{\ast}\|\le\|p_{n}bp^{\perp}_{n}\|+\|p^{\perp}_{n}bp_{n}\|,$$	
we have $\|P_n(b)-P_n(b^{\ast})^{\ast}\|\to 0$ as $n\to\infty$.

(a) $\Leftarrow$ (b) Assume that (b) holds. For any $\epsilon>0$,  there is a $N$ such that for any $n>N$ we have
$$\|P_{n}(d)-P_{n}(d^{\ast})^{\ast}\|<\epsilon.$$
Then we find
\begin{equation*}
    \begin{aligned}
        \|p_{n}^{\perp}dp_{n}\|&=\left\|\begin{bmatrix}
        	0       &  0\\
        	-d_{21}  &  0\\
          \end{bmatrix}\right\|  \le \left\|\begin{bmatrix}
        	0       &  d_{12}\\
        	-d_{21}  &  0\\
        \end{bmatrix}\right\|\\
        &= \left\|P_{n}(d)-P_{n}(d^{\ast})^{\ast}\right\|\le \epsilon.
      \end{aligned}
\end{equation*}
Hence  $\|p_{n}^{\perp}dp_{n}\|\to 0$ as $n\to\infty$. We can prove that  $\|p_{n}dp_{n}^{\perp}\|\to 0$ similarly.
Finally we have $\|[d,p_{n}]\|\to0$ as $n\to\infty$.
\qed

Such a sequence $\{P_n\}$ of Rota-Baxter operators of weight $-1$ on $B(\mathcal{H})$ in (b) of the above proposition is called a quasi symmetric Rota-Baxter operator sequence on $B(\mathcal{H})$.



\end{document}